\documentclass[10pt,a4paper,twoside]{article}
\usepackage{epsfig}
\usepackage[latin1]{inputenc}
\usepackage{amsmath}
\usepackage{amsfonts}
\usepackage{amssymb}

\newcommand{\N}{\mathbb N}  
\newcommand{\PP}{\mathbb P}

\newcommand{\G}{\mathbb G}
\newcommand{\A}{\mathbb A}

\newcommand{\vtl}{\, | \,}

\newcommand{\join}{\vee}

\newcommand{\lra}{\longrightarrow}

\newcommand{\dra}{\dashrightarrow}

\newcommand{\al}{\alpha}

\newcommand{\si}{\sigma}

\newcommand{\De}{\Delta}

\newcommand{\Si}{\Sigma}
\newcommand{\Tht}{\Theta}

\newtheorem{coro}{Corollary}
\newtheorem{defin}{Definition}

\newtheorem{lemma}{Lemma}

\newtheorem{prop}{Proposition}

\newtheorem{theo}{Theorem}

\newenvironment{proof}{{\em Proof:}}{\hfill\rule{2mm}{2mm}}

\begin{document}

\title{Trisecant Lemma for Non Equidimensional Varieties\\}
\author{Jeremy Yirmeyahu Kaminski$^1$, Alexei Kanel-Belov$^2$ and Mina Teicher$^3$ \\ \\
$^1$ Holon Institute of Technology, \\
Department of Computer Science, \\
Holon, Israel. \\
E-mail: kaminsj@hait.ac.il \\ \\
$^2$ The Hebrew University, \\
Department of Mathematics,  \\
Jerusalem, Israel. \\
E-mail: kanel@mccme.ru \\ \\
$^3$ Bar-Ilan University, \\
Department of Mathematics, \\
Ramat-Gan, Israel. \\
E-mail: teicher@macs.biu.ac.il
}

\date{}
\maketitle

\begin{abstract}
The classic trisecant lemma states that if $X$ is an integral curve of $\PP^3$ then the variety of trisecants has dimension one, unless the curve is planar and has degree at least $3$, in which case the variety of trisecants has dimension 2. Several generalizations of this lemma has been considered \cite{Laudal-77, Ran-91, Adlandsvik-87, Adlandsvik-88, Flenner-all-99}. In \cite{Laudal-77}, the case of an integral curve embedded in $\PP^3$ is further investigated leading to a result on the planar sections of the such a curve. On the other hand, in \cite{Ran-91}, the case of higher dimensional varieties, possibly reducible, is inquired. For our concern, the main result of \cite{Ran-91} is that if $m$ is the dimension of the variety, then the union of a family of $m+2$-secant lines has dimension at most $m+1$. A further generalization of this result is given in \cite{Adlandsvik-87, Adlandsvik-88, Flenner-all-99}. In this latter case, the setting is the following. Let $X$ in an irreducible projective variety over an algebraically closed field of characteristics zero. For $r \geq 3$, if every $(r-2)-$plane $\overline{x_1 ... x_{r-1}}$, where the $x_i$ are generic points, also meets $X$ in an $r-$th point $x_r$ different from $x_1,...,x_{r-1}$, then $X$ is contained in a linear subspace $L$, with $codim_L X \leq r-2$.

In this paper, our purpose is first to present another derivation of this result for $r=3$ and then to introduce a generalization to non-equidimensional varities. For the sake of clarity, we shall reformulate our first problem as follows. Let $Z$ be an equidimensional variety (maybe singular and/or reducible) of dimension $n$, other than a linear space, embedded into $\PP^r$, $r \geq n+1$. The variety of trisecant lines of $Z$, say $V_{1,3}(Z)$, has dimension strictly less than $2n$, unless $Z$ is included in a $(n+1)-$dimensional linear space and has degree at least $3$, in which case $\dim(V_{1,3}(Z)) = 2n$. This also implies that if $\dim(V_{1,3}(Z)) = 2n$ then $Z$ can be embedded in $\PP^{n+1}$.

Then we inquire the more general case, where $Z$ is not required to be equidimensional. In that case, let $Z$ be a possibly singular variety of dimension $n$, that may be neither irreducible nor equidimensional, embedded into $\PP^r$, where $r \geq n+1$, and $Y$ a proper subvariety of dimension $k \geq 1$. Consider now $S$ being a component of maximal dimension of the closure of $\{l \in \G(1,r) \vtl \exists p \in Y, q_1, q_2 \in Z \backslash Y, q_1,q_2,p \in l\}$. We show that $S$ has dimension strictly less than $n+k$, unless the union of lines in $S$ has dimension $n+1$, in which case $dim(S) = n+k$. In the latter case, if the dimension of the space is stricly greater then $n+1$, the union of lines in $S$ cannot cover the whole space. This is the main result of our work. We also introduce some examples showing than our bound is strict.
\end{abstract}


%=========================================
%=========================================
\section{Introduction}

The classic trisecant lemma states that if $X$ is an integral curve of $\PP^3$ then the variety of trisecants has dimension one, unless the curve is planar and has degree at least $3$, in which case the variety of trisecants has dimension 2. Several generalizations of this lemma has been considered \cite{Laudal-77, Ran-91, Adlandsvik-87, Adlandsvik-88, Flenner-all-99}. In \cite{Laudal-77}, the case of an integral curve embedded in $\PP^3$ is further investigated leading to a result on the planar sections of the such a curve. On the other hand, in \cite{Ran-91}, the case of higher dimensional varieties, possibly reducible, is inquired. For our concern, the main result of \cite{Ran-91} is that if $m$ is the dimension of the variety, then the union of a family of $m+2$-secant lines has dimension at most $m+1$. A further generalization of this result is given in \cite{Adlandsvik-87, Adlandsvik-88, Flenner-all-99}. In this latter case, the setting is the following. Let $X$ in an irreducible projective variety over an algebraically closed field of characteristics zero. For $r \geq 3$, if every $(r-2)-$plane $\overline{x_1 ... x_{r-1}}$, where the $x_i$ are generic points, also meets $X$ in an $r-$th point $x_r$ different from $x_1,...,x_{r-1}$, then $X$ is contained in a linear subspace $L$, with $codim_L X \leq r-2$.

In this paper, our purpose is first to present another derivation of this result for $r=3$ and then to introduce a generalization to non-equidimensional varities. For the sake of clarity, we shall reformulate our first problem as follows. Let $Z$ be an equidimensional variety (maybe singular and/or reducible) of dimension $n$, other than a linear space, embedded into $\PP^r$, $r \geq n+1$. The variety of trisecant lines of $Z$, say $V_{1,3}(Z)$, has dimension strictly less than $2n$, unless $Z$ is included in a $(n+1)-$dimensional linear space and has degree at least $3$, in which case $\dim(V_{1,3}(Z)) = 2n$. This also implies that if $\dim(V_{1,3}(Z)) = 2n$ then $Z$ can be embedded in $\PP^{n+1}$.

Then we inquire the more general case, where $Z$ is not required to be equidimensional. In that case, let $Z$ be a possibly singular variety of dimension $n$, that may be neither irreducible nor equidimensional, embedded into $\PP^r$, where $r \geq n+1$, and $Y$ a proper subvariety of dimension $k \geq 1$. Consider now $S$ being a component of maximal dimension of the closure of $\{l \in \G(1,r) \vtl \exists p \in Y, q_1, q_2 \in Z \backslash Y, q_1,q_2,p \in l\}$. We show that $S$ has dimension strictly less than $n+k$, unless the union of lines in $S$ has dimension $n+1$, in which case $dim(S) = n+k$. In the latter case, if the dimension of the space is stricly greater then $n+1$, the union of lines in $S$ cannot cover the whole space. This is the main result of our work. We also introduce some examples showing than our bound is strict. 

The methods we use to prove these results are purely algebraic and are valid over any algebraically closed field of characteristic zero. Our reasoning consists basically in inquiring first the local implications on the tangent spaces for the trisecant lines variety being of full dimension. Then the global result is deduced using the so-called Terracini's lemma~\cite{Zak-93}.

The paper is organized as follows. We should first recall some standard material in order to fix terminology and notations, in section~\ref{sec::background}. Then we come to our results, in section~\ref{sec::core}. More precisely, in subsection~\ref{subsec::equidim}, the case of equidimensional varieties is investigated, while in subsection~\ref{subsec::gencase}, we deal with the more general case.

%=========================================
%=========================================
\section{Notations and Background}
\label{sec::background}

In this section, we recall some standard material on incident varieties, that will be used in the sequel.

%=========================================
\subsection{Variety of Incident Lines}

Let $\G(1,n) = G(2,n+1)$ be the Grassmannian of lines included in $\PP^n$. Remind that $\G(1,n)$ can be canonically embedded in  $\PP^{N_1}$, where $N_1=\binom{2}{n+1}-1$, by the Pl\"ucker embedding and that $\dim(\G(1,n))=2n-2$. Hence a line in $\PP^n$ can be regarded as a point in $\PP^{N_1}$, satisfying the so-called Pl\"ucker relations. These relations are quadratic equations that generate a homogeneous ideal, say $I_{\G(1,n)}$, defining $\G(1,n)$ as a closed subvariety of $\PP^{N_1}$. Similarly the Grassmannian, $\G(k,n)$, gives a parametrization of the $k-$dimensional linear subspaces of $\PP^n$. As for $\G(1,n)$, the Grassmannian $\G(k,n)$ can be embedded into the projective space $\PP^{N_k}$, where $N_k=\binom{k+1}{n+1}-1$. Therefore for a $k-$dimensional linear subspace, $K$, of $\PP^n$, we shall write $[K]$ for the corresponding projective point in $\PP^{N_k}$. 

\begin{defin}
Let $X \subset \PP^n$ be an irreducible variety. We define the following variety of incident lines:
$$
\De(X) = \{l \in \G(1,n) \vtl l \cap X \neq \emptyset \}.
$$
\end{defin}

The codimension $c$ of $X$ and the dimension of $\De(X)$ are related by the following lemma.

\begin{lemma}
Let $X \subset \PP^n$ be an irreducible closed variety of codimension $c \geq 2$. Then $\De(X)$ is an irreducible variety of $\G(1,n)$ of dimension $2n-1-c$.
\end{lemma}
\begin{proof}
Consider the following incidence variety $\Si = \{(l,p) \in \G(1,n) \times X \vtl p \in l\} \subset \De(X) \times X$, endowed with the canonical projections $\pi_1: \Si \lra \De(X)$ and $\pi_2: \Si \lra X$. The generic fiber of $\pi_1$ is finite (otherwise it is clear that $X = \PP^n$). Thus $\dim(\Si) = \dim(\De(X))$. For all $p \in X$, the fiber $\pi_2^{-1}(p)$ is isomorphic to $\PP^{n-1}$ and has dimension $n-1$. Therefore $\Si$ is irreducible and has dimension $n-c+n-1 = 2n-c-1$ as shown in \cite{Harris-92,Shafarevich-77}. Since $\pi_1$ is surjective and continuous (in Zariski topology), $\De(X)$ is also irreducible and has dimension $2n-c-1$.
\end{proof}

The following simple result will be useful in the sequel.

\begin{lemma}
\label{lemma::no_delta_inclusion}
Let $X_1$ and $X_2$ be two irreducible closed varieties in $\PP^n$ of codimension greater or equal to $2$. Then $\De(X_1) \not\subset \De(X_2)$ unless $X_1 \subset X_2$. 
\end{lemma}
\begin{proof}
Assume that $\De(X_1) \subset \De(X_2)$ and $X_1 \not\subset X_2$. Consider a point $p \in X_1 \backslash X_2$. and an hyperplane $H$ not passing through $p$. Consider the projection $\pi: \PP^n \backslash \{p\} \lra H, q \mapsto \overline{qp} \cap H$, which maps a point $q \in \PP^n \backslash \{p\}$ to the point of intersection of the line $\overline{qp}$ with the hyperplane $H$. The projection is surjective and so is $\pi_{|X_2}$, because $\De(X_1) \subset \De(X_2)$. Thus $\dim(X_2) \geq n-1$, which is impossible, because $co\dim(X_i) \geq 2$ for each $i$.
\end{proof}

%=========================================
\subsection{Join Varieties}

Consider $m < n$ closed irreducible varieties $\{Y_i\}_{i=1,...,m}$ embedded in $\PP^n$, with codimensions $c_i \geq 2$. Consider the {\it join variety}, $J = J(Y_1,...,Y_m) = \De(Y_1) \cap ... \cap \De(Y_m)$, included in $\G(1,n)$. We assume that $\sum_{i=1,...,m} c_i \leq 2n-2+m$, so that $J$ is not empty. We shall first determine the irreducible components of $J$. 

Let $U$ be the open set of $Y_1 \times ... \times Y_m$ defined by $\{(p_1,...,p_m) \in Y_1 \times ... \times Y_m \vtl \exists i \neq j, p_i \neq p_j \}$.  Let $V$ be the locally closed set made of the $m-$tuples in $U$, which points are collinear. Let $s: V \lra \G(1,n)$ be the morphism that maps a $m-$tuple of aligned points to the line they generate. Let $S \subset \G(1,n)$ be the closure of the image of $s$. 

First let us look at the irreducible components of $S$. These components could be classified in several classes according to the number of distinct points in the $m-$tuples that generate them. For example consider the case where $m=3$. The locally closed subset of $Y_1 \times Y_2 \times Y_3$, made of triplets of three distinct and collinear points, generates one component of $S$. Now, if $Y_{12}$ is an irreducible component of $Y_1 \cap Y_2$ not contained in $Y_3$, then the lines generated by a point of $Y_{12} \backslash Y_3$ and another point in $Y_3$ form also an irreducible component of $S$. Also let $Z$ be an irreducible component of $Y_1 \cap Y_2 \cap Y_3$, then the lines generates by a point of $Z$ and another point in $Y_1$ are the intersection of the secant variety of $Y_1$ with $\De(Z)$, and form an irreducible component of $S$ too. In the general case, the following lemma will be enough for our purpose. 

\begin{lemma}
The irreducible components of $J$ are:
\begin{enumerate}
\item $\De(Z)$, where $Z$ runs over all irreducible components of $Y_1 \cap ... \cap Y_m$,
\item the irreducible components of $S$, which are not included in any component of the form $\De(Z)$.
\end{enumerate}
\end{lemma}
\begin{proof}
These sets are all irreducible closed subsets of $J$. There is a finite number of such sets and their union covers $J$. Thus the irreducible components of $J$ are certainly some of these sets. 

For every irreducible component, $Z$, of $Y_1 \cap ... \cap Y_m$, $\De(Z) \not\subset S$ (otherwise by proceeding similarly than in lemma~\ref{lemma::no_delta_inclusion}, we have $\dim(Y_1 \cup ... \cup Y_m) \geq n-1$). Now by lemma~\ref{lemma::no_delta_inclusion}, $\De(Z_1) \not\subset \De(Z_2)$ for any two irreducible components, $Z_1$ and $Z_2$, of $Y_1 \cap ... \cap Y_m$. Therefore $\De(Z)$ is a maximal irreducible closed subset of $J$ for every irreducible component $Z$ of $Y_1 \cap ... \cap Y_m$. 

Every irreducible component $S_1$ of $S$, not included in any component of the form $\De(Z)$, is also a maximal irreducible closed subset of $J$. 
\end{proof}

For simplicity, we shall call the irreducible components of $S$ {\it joining components} of $J$ and components of the form $\De(Z)$ for some irreducible component $Z$ of $Y_1 \cap ... \cap Y_m$, {\it intersection components}. \\

We conclude this section by quoting Terracini's lemma, in the form we shall use it later. For this purpose and throughout the paper, we use the following notations. If $X$ is a projective subvariety of $\PP^n$, we shall write $T_p(X)$ for the projective embedded tangent space of $X$ at $p$. The Zariski tangent space is denoted $\Tht_p(X)$. Let $CX$ be the affine cone over $X$, then $T_p(X)$ is the projective space of one-dimensional subspaces of $\Tht_q(CX)$, where $q \in \A^{n+1}$ is any point lying over $p$. Hence for a morphism $f$ between two projective varieties $X$ and $Y$, which can be also be viewed as a morphism between $CX$ and $CY$, the differential  $df_p: T_p(X)\backslash \PP(\ker(\phi)) \lra T_{f(p)}(Y)$ is induced by the differential $\phi$ between the Zariski tangent spaces $ \Tht_q(CX)$ and $\Tht_{f(q)}(CY)$. For simplicity, we shall write: $df_p: T_p(X) \lra T_{f(p)}(Y)$, while it is understood that $df_p$ might be defined on a proper subset of $T_p(X)$.

\begin{lemma} {\bf Terracini's Lemma.} 
\label{lemma::Terracini}
Let $X$ and $Y$ be two irreducible projective varieties embedded in $\PP^n$, over an algebraically closed field of characteristic zero. Let $W(X,Y)$ be the union of the lines in $J(X,Y)$. Let $z$ be a point in $W(X,Y)$ lying neither on $X$ nor on $Y$. Then the tangent space of $W(X,Y)$ at $z$ is given by the following equality:
$$
T_z(W(X,Y)) = \langle T_x(X), T_y(Y) \rangle,
$$
where $(x,y) \in X \times Y$, such that $z \in \langle x,y \rangle = \overline{xy}$ and  $\langle \mbox{ } \rangle$ denotes the linear span.
\end{lemma}

A slightly more general statement and a proof can be found in \cite{Zak-93}.

%=========================================
%=========================================
\section{Generalizations Of The Trisecant Lemma}
\label{sec::core}

In this section, we shall introduce two generalization of the trisecant lemma. The first one is about equidimensional varieties, while the second one deals with a more general situation. 

In terms of join varieties, the classical trisecants lemma and the generalizations we introduce, are related to join components. A similar treatment of intersection component is easy to give and is summarized in the following lemma, immediatly deduced from lemma~\ref{lemma::no_delta_inclusion}.

\begin{lemma}
Let $Y_1$ and $Y_2$ be two distinct irreducible varieties embedded in some projective space. Let $Y$ be a third irreducible variety. $\De(Y)$ cannot contain any intersection component $\De(Z)$ of $J(Y_1,Y_2) = \De(Y_1) \cap \De(Y_2)$, unless $Z \subset Y$.
\end{lemma}

Before we proceed, we shall prove some results, useful in the sequel.

%======================
\subsection{Preliminary Properties}

The following proposition also illustrates the techniques we use in the paper. It can be viewed as a generalization of a well-known result of Samuel,~\cite{Hartshorne-77} page 312, which deals with smooth curves.

\begin{prop}
\label{prop::pencil}
Let $X$ be an irreducible closed subvariety of $\PP^n$ of dimension $k$. If there exists $L \in \G(k-1,n)$, such that for all  points $p \in U_0$, where $U_0$ is a dense open set of $X$, $L \subset T_p(X)$, then $X$ is a $k-$dimensional linear space, containing $L$.
\end{prop}
\begin{proof}
Let ${\mathcal T}X$ be the closure of $\{[T_p(X)] \vtl p \in X, p \mbox{ regular}\}$ in $\G(k,n)$. ${\mathcal T}X$ is the closure of the image of a dense open set of $X$ by the Gauss map. Therefore ${\mathcal T}X$ is irreducible. Consider the following rational map $X \dra \G(k,n), p \join L$, where $\join$ is the join operator \cite{Barnabei-all-85}, equivalent to the classical exterior product \footnote{As in \cite{Barnabei-all-85}, the departure from the classical notation is amply justified by the geometric meaning on the operator.}. Let $\si_L$ be the subvariety of $\G(k,n)$ made of the linear spaces that contains $L$. Thus $\dim(\si_L) = n-k$. 

Let $U$ be the open set of $X$ made of the regular points of $U_0$ which do not lie on $L$. Consider the morphism: $f: U \lra \si_L, p \mapsto p \join L$. For each $p \in U$,  $f(p)$ is simply the tangent space of $X$ at $p$. Therefore the image of $f$ is dense in ${\mathcal T}X$. 

Since the ground field is assumed to have characteristic zero, there exists a dense open set $V$ of $X$ such that for any point $p$ in $V$, the differential $df_p$ is surjective,~\cite{Hartshorne-77} page 271. 

This differential is simply: $df_p: T_p(X) \lra T_{f(p)}({\mathcal T}X), a \mapsto a \join L$. Therefore $df_p$ is constant over $T_p(X) \backslash L$ and takes the value $[T_p(X)] = df_p(p)$. Thus $\dim({\mathcal T}X) = 0$. Since ${\mathcal T}X$ is irreducible, it is a single point corresponding to a $k-$dimensional linear space, say $T$, containing $L$. Finally $X \subset T$, $\dim(X) = k$ and $X$ is closed, therefore we have: $X=T$.
\end{proof}

Note that this fact does not hold in positive characteristic as the following example shows. Consider the curve in $\PP^3$, over a field $K$ of characteristic $p$, defined by the ideal $<\!y^p-zt^{p-1},x^p-yt^{p-1}\!> \subset K[x,y,z,t]$, with $t=0$ being the plane at infinity. The tangent space at $(x_0,y_0,z_0,t_0)$ is given by the following system of linear equations: $\{t_0^{p-1}z+(p-1)z_0t_0^{p-2}t=0, t_0^{p-1}y+(p-1)y_0t_0^{p-2}t=0\}$. Every two tangent spaces are parallel and therefore they all contain the same point at infinity. However the the curve is not a line. Note that the point $(0,0,1,0)$ is a singular point of the curve.

The next proposition is used throughout the paper several times. The underlying idea is the following. Let $L$ be a $k-$dimensional linear space. If the tangent space to an irreducible variety at a generic point always spans with $L$ a $(k+1)-$dimensional linear space, then the variety itself must be included into a $(k+1)-$dimensional linear space, containing $L$.

\begin{prop}
\label{prop::planar}
Let $X$ be an irreducible closed subset of $\PP^n$, with $\dim(X) = r$. If there exists $L \in \G(k,n)$, such that for all points $p \in U_0$, where $U_0$ is dense open set of $X$, $\dim(L \cap T_p(X)) \geq r-1 $, then $X$ is included in a $(k+1)-$dimensional linear space, containing $L$.
\end{prop}
\begin{proof}
If $X \subset L$, then there is nothing to prove. Therefore let us assume that $X \not\subset L$. Let $\si_L \subset \G(k+1,n)$ be the set of $(k+1)-$dimensional linear spaces that contains $L$. Consider the rational map: $f: X \dra \si_L, p \mapsto p \join L$. This mapping is defined over the open set $U$ of regular points in $(X \backslash L) \cap U_0$. Each such point is mapped to the $(k+1)-$dimensional space generated by $p$ and $L$. Since $\dim(T_p(X) \cap L) = r-1$, we have the following inclusion $T_p(X) \subset p \join L = f(p)$, for $p \in U$. 
Let $Y$ be the closure of $f(U)$ in $\si_L$. Thus $Y$ is irreducible. 

Since the ground field is assumed to have characteristic zero, there exists a dense open set $V$ of $X$ such that for any point $p$ in $V$, the differential $df_p$ is surjective,~\cite{Hartshorne-77} page 271. 

This differential is simply: $df_p: T_p(X) \lra T_{f(p)}(Y), a \mapsto a \join L$. Since $T_p(X) \subset p \join L$, $df_p$ is constant over $T_p(X) \backslash L$ and takes the value $p \join L = df_p(p)$. Thus $\dim(Y) = 0$. Since $Y$ is irreducible, $Y$ is a single point corresponding to a $(k+1)-$dimensional linear space, say $K$, containing $L$. Therefore $X \subset K$.
\end{proof}

This proposition does not hold in positive characteristic. Indeed over a field of characteristic $p$, for the curve in $\PP^3$ defined by the following ideal:\\ $<\!yt^{p-1}-x^p,zt^{p^2-1}-x^{p^2}\!>$, all the tangent lines are parallel and therefore intersect in some point at infinity. But the curve is not a line.

Before we come to investigate our initial question, let us first show in the case of two varieties embedded in $\PP^n$ with $n \geq 3$, that the join has necessarily a unique joining component, which has the expected dimension.

\begin{lemma}
\label{lemma::maxdim}
Let $Y_1$ and $Y_2$ be two distinct irreducible varieties embedded in $\PP^n$. Let $c_i \geq 2$ be the codimension of $Y_i$. Then the join $J=J(Y_1,Y_2)$ has a unique joining component $S$, which dimension is $2n-(c_1+c_2)$.
%, which is the dimension of $J$ itself. 
\end{lemma}
\begin{proof}
%The join $J$ has obviously dimension $2n-(c_1+c_2) = 2n-1-c_1 + 2n-1-c_2 - (2n-2)$. 
Let $\De = \{(y_1,y_2) \in Y_1 \times Y_2 \vtl y_1 = y_2\}$. Let $U$ be the open set of $Y_1 \times Y_2$, defined as $U = (Y_1 \times Y_2) \backslash \De$. Let $s: U \lra \G(1,n), (p,q) \mapsto \overline{pq}$ be the morphism, which maps a couple of points in $U$ to the line they generate. Let $S$ be the closure of $s(U)$ in $\G(1,n)$. As $U$ is irreducible, so is $S$. It is therefore the unique joining component of $J$. The general fiber is finite. Thus $\dim(S) = \dim(U) = \dim(Y_1 \times Y_2) = 2n-(c_1+c_2)$. 
\end{proof}

Eventually, we also the following lemma, which be useful in the sequel.

\begin{lemma}
\label{lemma::fibers}
Let $Y_1$ and $Y_2$ be two irreducible varieties embedded in $\PP^n$, with dimensions $d_1,d_2$ both smaller or equal to $n-2$. Let $S$ be the unique joining component of $J=J(Y_1,Y_2)$. Then $dim(S) = s = d_1+d_2$.

The union of the lines in $S$ is an irreducible variety of dimension strictly greater than $max(d_1,d_2)$. For a generic point $p$ in $Y_i$, the dimension of the variety of lines in $S$ passing through $p$ is $d_{3-i}$.

Moreover, if there exists an irreducible variety $Y$ of dimension $d \leq max(d_1,d_2)$, such that $S \subset \De(Y)$, then $d=max(d_1,d_2)$ and for a generic point $p$ in $Y$, the dimension of the lines in $S$ passing through $p$ is $min(d_1,d_2)$.  
\end{lemma}
\begin{proof}
We shall assume, without loss of generality, that $d_1 \geq d_2$.

{\it Step 1 - } Consider first the following incidence variety:
$$
\Si = \{(l,p) \in S \times \PP^n \vtl p \in l \},
$$
endowed with the canonical projections: $\pi_1: \Si \lra S$ and $\pi_2: \Si \lra \PP^n$. For all $l \in S$, the fiber $\pi_1^{-1}(l)$ is irreducible and has dimension $1$ and $S$ is irreducible. Thus $\Si$ is irreducible and
$\dim(\Si) = \dim(S) + 1 = s+1$. Let $W = \pi_2(\pi_1^{-1}(S)) = \bigcup_{l \in S} l$. Then $W$ is irreducible, since $\Si$ is irreducible. Since for each $i$, $Y_i \subset W$, $\dim(W) \geq max(d_1,d_2)$. Furthermore the generic fiber of $\pi_2$ has dimension less or equal to $d_2$; indeed, the fiber at a generic point $p$ is included in $\{(\overline{qp},p) \vtl q \in Y_2\}$. Thus $\dim(W) > max(d_1,d_2)$. \\

{\it Step 2 - } For $p \in Y_1$, consider the open set $U=Y_2 \backslash \{p\}$ and the morphism $f: U \lra S, q \mapsto \overline{pq}$. Since $U$ is irreducible, then so is $f(U)$. For a generic line $l$ in $f(U)$, the fiber $f^{-1}(l)$ is finite (otherwise $Y_2$ is a cone with vertex $p$, which is impossible for a generic $p \in Y_1$). Therefore $\dim(f(U)) = dim(Y_2)$. A similar conclusion is valid for a generic point of $Y_2$. Therefore the dimension of the lines in $S$ passing through a general point in $Y_i$ is $d_{3-i}$.\\

{\it Step 3 - } For a point $p \in Y$, let $X_p$ be the variety of lines in $S$ passing through $p$. Let $Z$ be the subvariety of $Y$, defined as the set of points for which $X_p$ is not empty. Thus $S \subset \De(Z)$. 

Let us show that $Z$ is irreducible. If $Z = E \cup F$, with $E$ and $F$ being closed subsets of $Z$, consider $S_1$ and $S_2$ being respectively the unique joining components of $J(E,Y_2)$ and $J(F,Y_2)$. Then $\dim(S_1) = \dim(E) + d_2$ and $\dim(S_2) = \dim(F) + d_2$. Moreover $S \subset S_1 \cup S_2$, therefore $max(\dim(E)+d_2,\dim(F)+d_2)) \geq d_1+d_2$, so that $max(\dim(E),\dim(F)) \geq d_1$. However $\dim(Z) \leq dim(Y) \leq d_1$. We conclude that either $E= Z$ or $F = Z$, and $\dim(Z) = d_1$. Thus $Z=Y_1$ and $\dim(Y) = d_1$. 

Let $S'$ be the unique joining component of $J(Z,Y_2)$. Then we have: $S \subset S'$. But $\dim(S) = \dim(S')$ and both varieties are irreducible closed varieties. Thus $S=S'$. Thus by a similar argument than in step 2, we get that for a generic point $p$ in $Y$, the dimension of $X_p$ is $d_2 = min(d_1,d_2)$. 
\end{proof}

%==============================================================================
\subsection{Equidimensional Varieties}
\label{subsec::equidim}

We are in a position to present our derivation of the general the trisecant lemma valid for equidimensional varieties. We shall first consider the following situation. Let $Y_1$ and $Y_2$ be two irreducible varieties embedded in $\PP^n$, for some $n \in 2 \N + 1$. Assume that $\dim(Y_1) = \dim(Y_2) = k =\frac{n-1}{2}$. The join $J(Y_1,Y_2)$ has necessarily a joining component $S$ of dimension $n-1$ as shown in lemma~\ref{lemma::maxdim}. We will show that if a third irreducible variety $Y$, of the same dimension, is such that $S \subset \De(Y)$, then the three varieties lie on the same $(k+1)-$dimensional linear subspace. Then we generalize to equidimensional varieties.

%=========================================
\subsubsection{Two Varieties of Equal Dimensions In A Space Which Dimension Is Odd} 

\begin{theo}
\label{theo::gencurves}
Let $n$ be an odd number. Consider two distinct irreducible closed varieties $Y_1$ and $Y_2$ in $\PP^3$, each of dimension $k = \frac{n-1}{2}$. By lemma~\ref{lemma::maxdim}, consider the joining component $S$ of $J(Y_1,Y_2)$, having dimension $n-1$. If there exists a third irreducible variety $Y$ of dimension $k$, distinct from $Y_1$ and $Y_2$, such that $S \subset \De(Y)$, then the three varieties lie in the same $(k+1)-$dimensional linear space, equal to the union of the lines in $S$.
\end{theo}
\begin{proof}
{\it Step 1 - } Let $W = \bigcup_{l \in S} l$. By lemma~\ref{lemma::fibers}, $W$ has dimension strictly greater than $k$. Moreover the same lemma shows that the dimension of the variety of lines in $S$ passing through a generic point $p$ in $Y$ has dimension $k$. \\

{\it Step 2 - } Let $l_0$ be a generic line in $S$. Let $q_i = l_0 \cap Y_i$ and $p_0 = l_0 \cap Y$. Since $l_0$ is
generic, these points can be assumed to be regular and $p_0 \not\in Y_1 \cup Y_2$. 

Let $\si_{p_0} \subset \G(1,n)$ be the set of lines passing through $p_0$. In general
$X_{p_0} = \si_{p_0} \cap S$ has dimension equal to $k$.

Consider now the morphism $f: Y_1 \lra \si_{p_0}, a \mapsto a \join p_0$. It is clear that $X_{p_0} \subset f(Y_1)$. The general fiber of $f$ being finite, $\dim(X_{p_0}) =\dim(Y_1)$ and $f(Y_1)$ being irreducible,  thus we have even the following equality: $X_{p_0} = f(Y_1)$. Therefore $f$ can be regarded as a morphism from $Y_1$~to~$X_{p_0}$: $f:~Y_1~\lra~X_{p_0}, a \mapsto a \join p$. Here again the expression of the differential of $f$ at $q_1$ is simply given by: $df_{p_1}: T_{q_1}(Y_1) \lra T_{l_0}(X_{p_0}), a \join p_0$. The line $l_0$ being generic, we shall assume that $\dim(T_{l_0}(X_{p_0}) = \dim(X_{p_0}) = k$.

Consider now $H_0 = \bigcup_{l \in T_{l_0}(X_{p_0})} l$. This linear space has dimension $k+1$.
The expression of $df_{p_1}$ shows that $T_{q_1}(Y_1) \subset H_0$.
Similarly we can deduce that $T_{q_2}(Y_2) \subset H_0$. Therefore the following inequality holds:
$\dim(T_{q_1}(Y_1) \cap T_{q_2}(Y_2)) \geq k-1$.

By the same reasoning, there exists a dense open set $U$  of $Y_1$, such that for each $q \in U$,
$\dim(T_q(Y_1) \cap T_{q_2}(Y_2)) \geq k-1$. \\

{\it Step 3 - } If $Y_2$ is a linear space of dimension $k$, then by proposition~\ref{prop::planar},
$Y_1$ is contained in a $(k+1)-$dimensional linear space containing $Y_2$.
A similar conclusion can be done if $Y_1$ is a linear space. \\

{\it Step 4 -} Assume now that neither $Y_1$ nor $Y_2$ is a linear space. Applying the reasoning than in step 2 to
$X_{q_1}$ and $X_{q_2}$, being respectively the set of lines in $S$ passing through $q_1$ and $q_2$, we get the following facts:
\begin{itemize}
\item there exists an open set $U_1$ of $Y_1$, such that for all
$q \in U_1$, $\dim(T_{q}(Y_1) \cap T_{p_0}(Y)) \geq k-1$,

\item there exists an open set $U_2$ of $U_2$, such that for all
$q \in V_0$, $\dim(T_{q}(Y_2) \cap T_{p_0}(Y)) \geq k-1$.
\end{itemize}

When $k=1$ (this is the case for curves in $\PP^3$), these inequalities just mean that the intersection are not empty. Then by proposition~\ref{prop::planar}, each $Y_i$ lies a $(k+1)-$dimensional linear
space $Q_i$ containing $T_{p_0}(Y)$. These two linear spaces $Q_1$ and $Q_2$ are identical, since there are both
generated by a line of $S$, namely $l_0$, and $T_{p_0}(Y)$. Let $Q$ denote this linear space.

Then $W$ being the union of the lines in $S$ is included in $Q$. Thus $Y$ is also included in $Q$. Then every line in $Q$ intersects the three varieties $Y_1, Y_2$ and $Y$. Therefore the Fano varieties of lines in $Q$ is the unique joining component of $J(Y_1,Y_2)$. The union of these lines is exactly $Q$. 
\end{proof}

%=====================================
\subsubsection{Generalized Trisecant Lemma for Equidimensional Varieties}

Since the proof is still valid, if some or all of the varieties $Y_1,Y_2$ and $Y$ are identical, we get a generalization of the trisecant lemma. We shall use the following notation: for a variety $X$, $V_{1,3}(X)$ is the closure in $\G(1,n)$ of 
$$
\{l \in \G(1,n) \vtl \exists p,q,r \in X, p \neq q, p \neq r, q \neq r, p,q,r \in l\}
$$.

\begin{theo}
\label{GenTriLemma1}
{\bf First Generalization of The Trisecant Lemma}\\
Let $Z$ be a possibly singular equidimensional variety (maybe reducible or not) of dimension $n$, other than a linear space, embedded into $\PP^r$, $r \geq n+1$. The variety of trisecant lines of $Z$, that is $V_{1,3}(Z)$, has dimension strictly less than $2n$, unless $Z$ is included in a $(n+1)-$dimensional linear space and has degree at least $3$, in which case $\dim(V_{1,3}(Z)) = 2n$.
\end{theo}
\begin{proof}
Two cases must be considered. 

{\it Step 1 - } If $r<2n+1$, then we can embed $\PP^r$ into $\PP^{2n+1}$ by a projective equivalence, so that we are in the setting of theorem~\ref{theo::gencurves}. Then the corollary follows immediately.

{\it Step 2 - } In the case where $r \geq 2n+1$, let us define $s=r-2n-1 \geq 0$. We shall prove the result by induction over $s$. If $s=0$, it is the content of the theorem~\ref{theo::gencurves}. 

Now it is left to show that if the result holds for some $s$, then it also holds for $s+1$. Let $p$ be a generic point in $\PP^r$, where $r=2n+1+s+1$, and let $H$ be any hyperplane in $\PP^r$, not passing through $p$. Let $Z'$ be the projection of $Z$ over $H$ through $p$. We can canonically identify $H$ to $\PP^{2n+1+s}$. Since the projection is generic and $dim(Z) < r-1$, the general fiber of the projection $\pi: Z \lra H$ is empty. However over $\pi(Z)$, the general fiber is finite. Therefore the dimension of $V_{1,3}(Z')$ is also $2n$. Then by the induction assumption, $Z'$ is included within a linear space $L' \subset H$ of dimension $n+1$. 

Let $L$ be the space generated by $p$ and $L'$. Then $\dim(L) = n+2$ and $Z \subset L$. Since $n+2 < 2n+1$, for $n > 1$, we can use the first step of the proof to conclude. Note that for $n=1$, the result can be easily deduced from the classic trisecant lemma.
\end{proof}

This result can also be expressed in the following terms.

\begin{coro}
Let $Z$ be a variety of dimension $n$. If the variety of trisecant lines $V_{1,3}(Z)$ has dimension $2n$, then $Z$ can be embedded into $\PP^{n+1}$.
\end{coro}

%=============================================
\subsection{Non-Equidimensional Case}
\label{subsec::gencase}

In this section, we turn a more general case. Our purpose is to generalize theorem~\ref{GenTriLemma1} to the case where the variety $Z$ is not equidimensional. As we proceeded before, we shall first inquire what happens with two irreducible varieties of complementary dimension.

\subsubsection{A Two Varieties Statement}

Let $Y_1$ and $Y_2$ be two irreducible closed varieties embedded in $\PP^n$. Let us assume that $\dim(Y_1) = k$ and $\dim(Y_2) = n-1-k$, where $\frac{n-1}{2} \leq k \leq n-2$. The varieties $Y_1$ and $Y_2$ are assumed to be distinct. Let $Y$ be another irreducible variety of dimension at most $k$, distinct from $Y_1$ and $Y_2$. By lemma~\ref{lemma::maxdim}, let $S$ be the joining component of $J(Y_1,Y_2)$, which dimension is $n-1$. Let $W$ be the subvariety of $\PP^n$ being the union of the lines in $S$. This setting is used throughout section~\ref{subsec::gencase}. Our purpose is to show that $W$ has dimension $k+1$. 

\paragraph{The Dimension of $Y$ is $k$}

\begin{lemma}
\label{lemma::dimY=k}
Let $Y_1$, $Y_2$ and $Y$ be varieties defined as just above. If $S \subset \De(Y)$, then the dimension of $Y$ must be equal to $k$.
\end{lemma}
\begin{proof}
It is clear by lemma~\ref{lemma::fibers}. 
\end{proof}

We are now in a position to address the determination of the dimension of $W$. 

\paragraph{$W$ Has Dimension $k+1$}

\begin{lemma}
\label{Lemma::tangent_spaces}
Let $Y_1,Y_2$ and $Y$ be varieties as in lemma~\ref{lemma::dimY=k}. Let $q_1$ and $q_2$ be generic points respectively on $Y_1$ and $Y_2$. Let $p(q_1,q_2) = \overline{q_1 q_2} \cap Y$ be an intersection point between the line $\overline{q_1 q_2}$ and the variety $Y$. The points $q_1,q_2$ and $p(q_1,q_2)$ can be assumed to be regular.

Then the tangent spaces $T_{q_1}(Y_1)$, $T_{q_2}(Y_2)$ and $T_{p(q_1,q_2)}(Y)$ 
%and the line $\overline{q_1 q_2}$ 
lie in the same $(k+1)-$dimensional linear space.
\end{lemma}
\begin{proof}
{\it Step 1 - } The points $q_1,q_2$ and $p(q_1,q_2)$ can indeed be assumed to be regular, since the set of singular points of an algebraic variety is a proper closed subvariety \cite{Shafarevich-77}. 

Let us prove first that the line $\overline{q_1 q_2}$ and the tangent spaces $T_{q_1}(Y)$ and $T_{p(q_1,q_2)}(Y)$ lie in the same $(k+1)-$dimensional linear space.

Let $\si_{q_2} \subset \G(1,n)$ be the set of lines passing through $q_2$. In general
$X_{q_2}$ has dimension equal to $k$ (by lemma~\ref{lemma::fibers}.

Consider now the morphism $f: Y_1 \lra \si_{q_2}, a \mapsto a \join q_2$. 
For each $a \in Y_1$, the line $a \join q_2$ lies in $S$. Therefore $f$ can be regarded as a morphism from $Y_1$~to~$X_{q_2}$: $f:~Y_1~\lra~X_{q_2}, a \mapsto a \join q_2$.
Again the differential of $f$ at $q_1$ is given as follows: $df_{q_1}: T_{q_1}(Y_1) \lra T_{\overline{q_1 q_2}}(X_{q_2}), a \mapsto a \join q_2$.

Consider now $H_{q_1,q_2} = \bigcup_{l \in T_{\overline{q_1 q_2}}(X_{q_2})} l$. This linear space has dimension $k+1$.
The expression of $df_{q_1}$ shows that $T_{q_1}(Y_1) \subset H_{q_1,q_2}$. Thus $H_{q_1,q_2}$ is the $(k+1)-$dimensional linear space generated by $T_{q_1}(Y_1)$ and the line $\overline{q_1 q_2}$: $H_{q_1,q_2} = \langle T_{q_1}(Y_1) , \overline{q_1 q_2} \rangle$, where $\langle \mbox{ } \rangle$ denotes the linear span as in Terracini's lemma. Similarly, one can prove that $T_{p(q_1,q_2)}(Y) \subset H_{q_1,q_2}$. \\

{\it Step 2 - } Consider now $\si_{p(q_1,q_2)}$, simply denoted $\si_p$ below, the set of lines passing through $p(q_1,q_2)$. 

Let $X_p = \si_p \cap S$. Lemma~\ref{lemma::fibers} shows that $\dim(X_p) = n-k-1$. Let $g:Y_2 \lra \si_p$ be the morphism which sends a point $a \in Y_2$ to the line $a \join p$, where $p=p(q_1,q_2)$. Since $X_p \subset g(Y_2)$, the general fiber of $g$ is finite, $g(Y_2)$ is irreducible and $\dim(Y_2) = \dim(X_p)$, the image of $g$ is simply $X_p$. Thus we can consider the morphism $g: Y_2 \lra X_p, a \mapsto a \join p$. The differential of $g$ at $q_2$ gives rise to the following morphism: $dg_{q_2}: T_{q_2}(Y_2) \lra T_{\overline{q_1 q_2}}(X_{p})$ given by $a \mapsto a \join p$. 

Let $K_{q_1,q_2} = \bigcup_{l \in T_{\overline{q_1 q_2}}(X_{p})} l$ be the union of lines in $T_{\overline{q_1 q_2}}(X_{p})$. It has dimension $n-k$. The expression of $dg_{q_2}$ shows that $T_{q_2}(Y) \subset K_{q_1,q_2}$. 

Now let $Z_1$ be the subvariety of $Y_1$ defined as follows: $Z_1 = \{q \in Y_1 \vtl \overline{qp} \in S\}$. It can be viewed as the trace over $Y_1$ of $X_p$. Let $h$ be the following morphism $h: Z_1 \lra X_p, a \mapsto a \join p$. Computing the differential of $h$ at $q_1$, shows that $T_{q_1}(Z_1) \subset K_{q_1,q_2}$. 

As $\dim(T_{q_1}(Z_1)) \geq n-k-1$, and in general $\overline{q_1q_2} \not\subset T_{q_1}(Z_1)$, and $\dim(K_{q_1,q_2}) = n-k$, we have $K_{q_1,q_2} = \langle T_{q_1}(Z_1) , \overline{q_1 q_2} \rangle$. Since $T_{q_1}(Z_1) \subset T_{q_1}(Y_1)$, we have: $K_{q_1,q_2} \subset H_{q_1,q_2}$ and therefore $T_{q_2}(Y_2) \subset H_{q_1,q_2}$.

Thus $T_{q_1}(Y_1)$, $T_{q_2}(Y_2)$ and $T_{p(q_1,q_2)}(Y)$ indeed linearly span a $(k+1)-$dimensional linear space.
\end{proof}

It is now possible to conclude using Terracini's lemma.

\begin{theo}
\label{theo::->W:k+1}
Let $Y_1,Y_2$ and $Y$ be varieties as in lemma~\ref{lemma::dimY=k}, then $W$ must have dimension $k+1$. 
\end{theo}
\begin{proof}
Consider smooth points $q_1 \in Y_1$ and $q_2 \in Y_2$. According to lemma~\ref{Lemma::tangent_spaces}, the tangent spaces $T_{q_1}(Y_1)$ and $T_{q_2}(Y_2)$ linearly span, together with the line $\overline{q_1q_2}$ a $(k+1)-$dimensional linear space, that we shall denote $K_{q_1,q_2}$. 

According to Terracini's lemma (lemma~\ref{lemma::Terracini}), the tangent space of $W$ at $\al q_1+q_2$, for some $\al \neq 0$ lies in $K_{q_1,q_2}$. Thus $\dim(W) \leq k+1$. Lemma~\ref{lemma::fibers} implies that $\dim(W) > k$. Therefore, we have: $\dim(W) = k+1$.
\end{proof}

In particular, the theorem shows that if $W$ covers all the space then there is no variety $Y$, distinct from $Y_1$ and $Y_2$ that intersects every line in $S$.

\paragraph{Example}

We shall now proceed to show how one can construct varieties as in section~\ref{subsec::gencase}. For any $k$, such that $\frac{n-1}{2} < k \leq n-2$, we can build varieties $Y_1$, $Y_2$ and $Y$ satisfying the following conditions:
\begin{itemize}
\item $\dim(Y_1) = \dim(Y) = k$,
\item $\dim(Y_2) = n-1-k$,
\item $J(Y_1,Y_2)$ has a joining component, $S$, of dimension $n-1$.
\item $S \subset \De(Y)$.
\end{itemize}

For this purpose, let $d = k-(n-1-k) = 2k-n+1 > 0$. Let $m > d$ be a natural number. Let $Z_1$ be a $d-$dimensional irreducible in $\A^m$, not passing through the origin. Let $Z_2$ be the single point variety made of the origin of $\A^m$. Let $f: \A^m \lra \A^m, (a_1,...,a_m) \mapsto (a_1/2,...,a_m/2)$. Let $Z = f(Z_1)$. Consider now $\hat{Y_1} = Z_1 \times \A^s$, $\hat{Y_2} = Z_2 \times \A^s$ and $\hat{Y} = Z \times \A^s$. 

If we take $s = k-d = n-k-1$ and $m = n-s = k+1 > d$, we have the following conditions: $\dim(\hat{Y_1}) = \dim(\hat{Y}) = k$, $\dim(\hat{Y_2}) = n-k-1$ and $\hat{Y_1}, \hat{Y_2}, \hat{Y} \subset \A^n$. 

Now we define $Y_1, Y_2, Y$ to be the projective closures of $\hat{Y_1}, \hat{Y_2}, \hat{Y}$. Then by lemma~\ref{lemma::maxdim}, we know that $J(Y_1,Y_2)$ has a joining component $S$ of dimension $n-1$. Moreover by construction we have: $S \subset \De(Y)$ and $W = \cup_{l \in S} l$ has dimension $k+1$.

\subsubsection{A General Statement}

The proof being true even when $Y_2 \subset Y_1$ and $Y_1 = Y$, we get the following consequence, which can be regarded as a generalization of the trisecant lemma, as well.

\begin{theo}
{\bf Second Generalization of the Trisecant Lemma} \\
Let $Z$ be a possibly singular variety of dimension $n$, that may be neither irreducible nor equidimensional, embedded into $\PP^r$, where $r \geq n+1$. Let $Y$ be a proper subvariety of $Z$ of dimension $k \geq 1$. Let $S$ an irreducible component of maximal dimension of $V_{1,3}(Y,Z)$, where $V_{1,3}(Y,Z)$ is the closure of $\{l \in \G(1,r) \vtl \exists p \in Y, q_1,q_2 \in Z \backslash Y, q_1 \neq q_2, p,q_1,q_2 \in l\}$. Then $S$ have dimension strictly less than $n+k$ unless the union of lines in $S$ has dimension $n+1$, in which case it has dimension $n+k$.
\end{theo}
\begin{proof}
{\it Step 1 -} The dimension of $S$ is at most $n+k$, since $n+k$ is exactly the dimension of the join $J(Y,Z)$. \\

{\it Step 2 - } If $r<n+k+1$, then we can embed $\PP^r$ into $\PP^{n+k+1}$ by a projective equivalence. Then according to theorem~\ref{theo::->W:k+1}, if $\dim(S)=n+k$ the union of lines in $S$ has dimension $n+1$. \\

{\it Step 3 -} If $r \geq n+k+1$, then let $s$ be $s=r-(n+k+1)$. If $s=0$, the result holds by theorem~\ref{theo::->W:k+1}. Assume now that the result is true for some $s \in \N$, let us prove it for $s+1$. 

The dimension $r$ of the space can be expressed as follows: $r=s+1+n+k+1$. Let $p$ be a generic point in $\PP^{s+1+n+k+1}$ and $H$ an hyperplane not passing through $p$. Then let $Z'$ (respectively $Y'$) be the projection of $Z$ (respectively $Y$) over $H$ through $p$. Then $Z'$ is embedded into a projective space of dimension $s+n+k+1$. The general fiber of the projection $\pi:Z \lra Z'$ is finite.

Each line in $S$ is projected over a line of the closure $V_{1,3}(Y',Z')$ of $\{l \in \G(1,r-1) \vtl \exists p \in Y', q_1,q_2 \in Z' \backslash Y', q_1 \neq q_2, p,q_1,q_2 \in l\}$. Let $S' \subset V_{1,3}(Y',Z')$ defined as being those lines which are built by the projection of lines in $S$. Since the general fiber of $\pi$ is finite, then $dim(S') = dim(S)$.

Therefore if $\dim(S) = n+k$ then $\dim(S') = n+k$. In that case, since $\dim(J(Y',Z')) = n+k$, $S'$ must be an irreducible component of maximal dimension of $V_{1,3}(Y'Z') \subset J(Y',Z')$. Thus by the induction assumption, $W' = \cup_{l \in S'} l$ has dimension $n+1$ and so $\dim(W) = n+1$, because the general fiber of $\pi: W \lra W'$ is finite. 
\end{proof}

Note that if $r > n+1$ and $dim(S) = n + k$, the theorem implies that the union of lines in $S$ cannot cover the whole space.

\paragraph{Example}

We shall now conclude by giving an example of a $n-$dimesional variety with $k-$secant lines variety of dimension $2n-1$, for $k \geq 3$. This improves the well known construction, also reminded in \cite{Ran-91}, of $n-$dimensional varieties admitting a $(n+1)-$dimensional variety of $k-$secant lines.

Let $p \in \A^3$ be the origin and consider an irreducible curve $X_1 \subset \A^3$ not passing through $p$. For $m \in \N, m \geq 2$, let $X_m$ be $f_m(X_1)$, where $f_m(x,y,z) = (mx,my,mz)$. For each $m \geq 1$, we shall denote $Y_m = X_m \times \A^{n-1}$. For a given $k \geq 3$, we shall define $Z_k = \cup_{1 \leq m \leq k} Y_m$. Then $\dim(Z_k) = n$ and $Z_k$ admits a family of $k-$secant lines, which dimension is $2n-1$.

We can also find an irreducible variety $Z$ containing $Z_k$ and having dimension $n' = n+1$. For this variety, the family of lines has dimension: $2n'-3$.

\bibliographystyle{amsplain}

\end{document}